\documentclass[11pt]{article}

\usepackage{amsfonts, amssymb}
\usepackage{amsmath}
\usepackage{color}
\usepackage{indentfirst}
\usepackage{verbatim}   % for \begin{comment}
\usepackage{mathrsfs}   % for \mathscr
\usepackage{yhmath}     % for anti-diagonal three consecutive dots

\newtheorem{theorem}{Theorem}[section]

\newtheorem{prop-def}{Proposition-Definition}[section]

\newtheorem{defi}[theorem]{Definition}

\begin{document}

\title{The 1882 Conjugacy Classes of the First  \\
Basic Renner Monoid of Type $E_6$
\footnote{Project supported by national NSF of China (No 11171202).} }
\date{}
\author{Zhenheng Li, Zhuo Li}
\maketitle

\vspace{ -1.2cm}

\vskip 7mm

\begin{abstract} In this paper we find all conjugacy classes of the first basic Renner monoid of type $E_6$. They are listed in the tables of the paper.

\vspace{ 0.3cm}
\noindent
{\bf Keywords:} Conjugacy, Renner monoid, Weyl group.

\vspace{ 0.3cm}
\noindent {\bf Mathematics Subject Classification 2010:}
20M32, 20M18.

\end{abstract}

\def\J {{\cal J}}
\def\a {\alpha}
\def\b {\beta}
\def\e {\varepsilon}
\def\s {\sigma}
\def\v {\vec}
\def\js {(\cal J, \sigma)}

\def\F {{\cal F}}
\def\ua {\underline a}
\def\ut {\underline t}
\def\e {\varepsilon}

\def\n { {\bf n} }
\def\r { {\bf r} }
\def\ek {\eta_{K}}
\def\eJ {\eta_{J}}
\newcommand{\be }{ \begin {}\\ \end{} }

\def\l {\langle}
\def\ra {\rangle}

\baselineskip 18pt

\section{Introduction}

The Renner monoids associated with exceptional $\mathcal J$-irreducible algebraic monoids of type $E_6, E_7, E_8$ are
closely related to Lie theory, semigroup theory, and other developments in mathematics and physics whose potential significance is unclear at this moment. The main result of this paper is the list of all conjugacy classes of the first basic Renner monoid of type $E_6$ (see Definition \ref{BasicMonoid}). This list is a result of an extensive computation on a Dell desktop computer OPTIPLEX 745 with a total physical memory 2,048.00 MB (about 1 week of computing time in May, 2012) using the method described in \cite{LLC2}. The code of the computation is written in GAP, a computer software system for algebraic computations \cite{GAP}.

Renner monoids $R$ were introduced by Renner \cite{R1} in his study of the structure of reductive monoids in 1986.
They play the same role for reductive monoids as the Weyl group does for groups in Lie theory. Significantly, $R$ is a finite inverse monoid and its unit group is a Weyl group $W$. Rook monoids, symplectic rook monoids and orthogonal rook monoids are classical examples of Renner monoids \cite{LR, PU1, R3}.

Two elements $a, b\in R$ are conjugate if and only if there exists $w\in W$ such that $b =  w a w^{-1}$. It is well known that two elements in a rook monoid are conjugate if and only if their cycle-link types are the same \cite{GK, GM, Ku, Lip, M2}. There is a one-to-one correspondence between conjugacy classes in a symplectic rook monoid and symplectic partitions \cite{CLL}. Reference \cite{LLC2} gives a systematic algorithm to determine the conjugacy of two elements of $R$ in terms of the action of certain parabolic subgroups on the cosets of some parabolic subgroups in the Weyl group. For integrity of the paper, we will briefly describe this algorithm in the next section. We refer the reader who is interested in other conjugacies in semigroup theory to \cite{GK, GM, Ku, KM1, KM2, PU3, PU4}, and those who are interested in conjugacy classes in semisimple algebraic groups and finite groups of Lie type to \cite{C1, C2} and \cite{H1}.

\section{The algorithm}
To describe the algorithm in \cite{LLC2} and to state our main results, we need some notation and basic facts. A linear algebraic monoid $M$ over an algebraic closed field is an affine algebraic variety together with an associative morphism from $M\times M$ to $M$ and an identity element $1\in M$. An irreducible algebraic monoid is a linear algebraic monoid whose underlying affine variety is not the union of two proper closed nonempty subsets. The unit group of $M$, consisting of all invertible elements in $M$, is an algebraic group. An irreducible algebraic monoid is reductive if its unit group is a reductive group.

Let $M$ be a reductive algebraic monoid, $T\subseteq G$ a maximal torus
of the unit group $G$, $B\subseteq G$ a Borel subgroup with
$T\subseteq B$, $N$ the normalizer of $T$ in $G$, $\overline N$ the
Zariski cloure of $N$ in $M$. Then $\overline N$ is a unit regular
inverse monoid which normalizes $T$, and $R=\overline N/T$ is a monoid, called the Renner monoid of $M$. This monoid contains the Weyl group $W = N/T$ as its unit group. Besides,
\[
    M = \bigsqcup_{r\in R} BrB, \quad \text{disjoin union}
\]
and if $s$ is a simple reflection then $BsB \cdot BrB \subseteq BsrB \cup BrB.$

Denote by $E(\overline T) = \{e \in \overline T \mid e^2 = e\}$ the set of idempotents in $\overline T$. It carries a partial order defined by
\[
    e \leq f \Leftrightarrow fe = e = ef.
\]
The set
$
\Lambda = \{e\in E(\overline T)\mid Be=eBe\}
$
is called the cross section lattice of $R$. Indeed,
\[
  R = \bigsqcup\limits_{e \in \Lambda} WeW,  \quad\mbox{disjoint union.}
\]
\noindent

A reductive monoid $M$ with zero 0 is called $\J$-irreducible if $\Lambda\setminus\{0\}$ has a unique minimal element. Let $G_0$ be a simple algebraic group and $\rho: G_0\to GL(V)$ be an irreducible rational representation over an algebraically closed field $K$. Let $K^* = K \setminus \{0\}$. Then
\[
M=\overline {K^*\rho(G_0)}
\]
is a $\J$-irreducible monoid, called the $\J$-irreducible monoid associated with $\rho$. Let $\mu_i$ ($1\le i \le l$) be the fundamental dominant weights of $G_0$ of type $X$, where $X = $ $A_l$, $B_l$, $C_l$, $D_l$, $E_6$, $E_7$, $E_8$, $F_4$ and $G_2$.

\begin{defi}\label{BasicMonoid}
The $\J$-irreducible monoid associated with $\mu_i$ is called the $i$-th basic monoid of type $X$; its Renner monoid is referred to as the $i$-th basic Renner monoid of type $X$.
\end{defi}

Denote by $\Delta$ the set of simple roots of $G$ relative to $T$, and let $S=\{s_\a \mid \a\in \Delta\}$ be the set of simple reflections that generate the Weyl group $W$. For each $e\in \Lambda$ let $W(e) = \{w\in W \mid we = ew\}$ and $W_*(e) = \{w\in W \mid we = ew = e\}$. Both $W(e)$ and $W_*(e)$ are parabolic subgroups of $W$, and $W_*(e)$ is a normal subgroup of $W(e)$. Let $W/W_*(e)$ be the set of left cosets of $W_*(e)$ in $W$, and let
\[
D_*(e)=\{w\in W \mid l(wu)=l(w)+l(u) \mbox{ for all } u\in W_*(e)\}.
\]
Then $D_*(e)$ is a set of left coset representatives of
$W/W_*(e)$, and each $w\in D_*(e)$ has a minimal length in $wW_*(e)$. Define the group action of $W(e)$ on $W/W_*(e)$ by
\[
    w\cdot uW_*(e) = wuw^{-1}W_*(e),
\]
where $w\in W(e)$ and $u\in W$. This action is well defined since $W_*(e)$ is a normal subgroup of $W(e)$.
We can now describe the algorithm for computing the conjugacy classes in a Renner monoid given in \cite{LLC2}.

\begin{theorem}\label{conjugaceTheoem}

{\rm a)} Each element in $R$ is conjugate to an element in $\{we \mid w\in D_*(e)\} \subseteq We$ for some $e\in \Lambda,$ and no element of $We$ is conjugate to an element of $Wf$ for $f\ne e.$

{\rm b)} If $e\in \Lambda$, then two elements $ue, ve$ in $We$ are conjugate if and only if the two cosets $uW_*(e)$ and $vW_*(e)$ are in the same $W(e)$-orbit in $W/W_*(e)$.
\end{theorem}

\section{Main results $E_6$}

With the help of GAP, we obtain that the first basic Renner monoid $R$ of type $E_6$ has 1882 conjugacy classes. The Dynkin diagram (c.f. \cite{H1}) in our calculation is as follows.
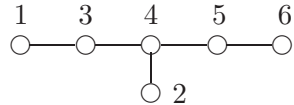
\begin{center}
\vskip 2mm
\setlength{\unitlength}{1mm}
\vskip 6mm
\setlength{\unitlength}{1mm}
\begin{picture}(80,20)
%\put(0, 18) {$E_6$:}
\put (24.8,19) {\circle{2.5}}
\put (33.5,19) {\circle{2.5}}
\put (42.1,19) {\circle{2.5}}
\put (51,19){\circle{2.5}}
\put (59.6,19) {\circle{2.5}}
\put (42.1,12.6) {\circle{2.5}}
\put (23.9,21.8) {$1$}
\put (32.5,21.8) {$3$}
\put (41.2,21.8) {$4$}
\put (50.3,21.8) {$5$}
\put (59,21.8) {$6$}
\put (45,11.3){$2$}
\put (26, 19.1) {\line (1,  0) {6}}
\put (34.7, 19.1) {\line (1,0) {6}}
\put (43.5, 19.1) {\line (1, 0) {6}}
\put (52.2, 19.1) {\line (1,0) {6}}
\put (42.1, 18) {\line (0, -1) {4}}
\end{picture}
\end{center}
\vskip -10mm
\noindent
The cross section lattice
$
\Lambda = \{0, ~e_0, ~e_1, ~e_2, ~e_3, ~e_4, ~e_5, ~e_6, ~e_7, ~e_8 = 1\}$
where $e_0$ is the minimal nonzero idempotent,
      $e_1 = \{s1\}$,
      ~$e_2 = \{s1, s3\}$,
      ~$e_3 = \{s1, s3, s4\}$,
      ~$e_4 = \{s1, s2, s3, s4\}$,
      ~$e_5 = \{s1, s2, s3, s4, s5\}$,
      ~$e_6 = \{s1, s3, s4, s5\}$,
      ~$e_7 = \{s1, s3, s4, s5, s6\}$,
      ~$e_8 = \{s1, s2, s3, s4, s5, s6\} = 1$.

\newpage
\paragraph{Description of tables.}
For simplicity we identify a conjugacy class with an individual representative element of that class, so we list only one representative for each class. Eleven tables are used to display necessary information. The first table is a summary of the classes associated with each idempotent; the rest display all the classes.

\vskip 5mm
\paragraph{Table $1$.} This table displays the number of conjugacy classes associated with each $e\in \Lambda$, and the parabolic subgroups $W^*(e)$ and $W_*(e)$. The notation $\l i, ..., j \ra$ in the table means the subgroup $\l s_i, ..., s_j\ra$ of $W$ generated by $s_i, ..., s_j$ where $1 \le i, j \le 6$.

{
\footnotesize
\begin{center}
\begin{tabular}{|c|c|c|c|}
\hline
\phantom{\rule{.2pt}{12pt}} % make more vertical space
$e$  & $W^*(e)$    & $W_*(e)$  &  Number of Classes \\
\hline\phantom{\rule{.4pt}{10pt}}
$0$         &  &  &  1 \\
\hline\phantom{\rule{.2pt}{10pt}}
$e_0$    & & $\langle 2, 3, 4, 5, 6 \rangle$ & $3$  \\
\hline\phantom{\rule{.2pt}{10pt}}
$e_1$    & $\langle 1 \rangle$ & $\langle 2, 4, 5, 6 \rangle$ & 16  \\
\hline\phantom{\rule{.2pt}{10pt}}
$e_2$    & $\langle 1, 3 \rangle$ & $\langle 2, 5, 6 \rangle$ & 113  \\
\hline\phantom{\rule{.2pt}{10pt}}
$e_3$    & $\langle 1, 3, 4 \rangle$ & $\langle 6 \rangle$ & 690  \\
\hline\phantom{\rule{.2pt}{10pt}}
$e_4$    & $\langle 1, 2, 3, 4 \rangle$ & $\langle 6 \rangle$ & 171  \\
\hline\phantom{\rule{.2pt}{10pt}}
$e_5$    & $\langle 1, 2, 3, 4, 5 \rangle$ & $\emptyset$ & 85  \\
\hline\phantom{\rule{.2pt}{10pt}}
$e_6$    & $\langle 1, 3, 4, 5 \rangle$ & $\emptyset$ & 628  \\
\hline\phantom{\rule{.2pt}{10pt}}
$e_7$    & $\langle 1, 3, 4, 5, 6 \rangle$ & $\emptyset$ & 150  \\
\hline\phantom{\rule{.2pt}{10pt}}
$e_8=1$    & $\langle 1, 2, 3, 4, 5, 6 \rangle$ & $\emptyset$ & 25  \\
\hline
\end{tabular}
%Number of conjugacy classes
\end{center}
}

\vskip 5mm
\paragraph{Table $i$.} The $i$-th table displays the representatives of classes associated with idempotent $e_{i-2}$ where $2\le i \le 10$. For convenience, the conjugacy class 0 is listed in Table 2. Notation $[ p, q, ..., r]$ in Table $i$ means
\[
    [p, q, ..., r] = s_ps_q ... s_r e_{i-2},  \text{ for } 1
\le p, \,q, \,..., \,r \le 6 \text{ and } i=2, ..., 10.
\]
A special case is $[ ~] = e_{i-2}$. For example, in Table 3,
\[
    [ 1, 3, 1, 4, 3 ] = s_1s_3s_1s_4s_3 e_1 \quad \text{and } \quad [~] = e_1.
\]

\vspace{5mm}
\paragraph{Table 2: Conjugacy Classes Associated with $0$ and $e_0$}
{
\footnotesize

\begin{center}
\begin{tabular}{|l|}
\hline
Conjugacy Classes \\
\hline\phantom{\rule{.2pt}{10pt}}
0	 \\
\hline\phantom{\rule{.2pt}{10pt}}
[ ]	 \\
\hline\phantom{\rule{.2pt}{10pt}}
[ 1 ]	 \\
\hline\phantom{\rule{.2pt}{10pt}}
[ 1, 2, 3, 4, 2, 5, 4, 3, 1, 6, 5, 4, 2 ] \\
\hline
\end{tabular}
\end{center}
}

\textheight 235mm
\topskip 0pt
\topmargin -10mm

\newpage
\paragraph{Table 3: Conjugacy Classes Associated with $e_1$}
{
\footnotesize
\begin{center}
\begin{tabular}{|l|}
\hline\phantom{\rule{.2pt}{10pt}}
Conjugacy Classes  \\
\hline\phantom{\rule{.2pt}{10pt}}
	[  ] \\
\hline\phantom{\rule{.2pt}{10pt}}
	[ 3 ] \\
\hline\phantom{\rule{.2pt}{10pt}}
	[ 1 ] \\
\hline\phantom{\rule{.2pt}{10pt}}
	[ 3, 1 ] \\
\hline\phantom{\rule{.2pt}{10pt}}
	[ 1, 3, 1, 4, 3 ] \\
\hline\phantom{\rule{.2pt}{10pt}}
	[ 2, 3, 1, 4, 2, 5, 4, 3, 6, 5 ] \\
\hline\phantom{\rule{.2pt}{10pt}}
	[ 3, 4, 2, 5, 4, 3, 6, 5, 4, 2 ] \\
\hline\phantom{\rule{.2pt}{10pt}}
	[ 2, 3, 1, 4, 2, 5, 4, 3, 1, 6, 5 ] \\
\hline\phantom{\rule{.2pt}{10pt}}
	[ 3, 4, 2, 5, 4, 3, 1, 6, 5, 4, 2 ] \\
\hline\phantom{\rule{.2pt}{10pt}}
	[ 3, 1, 4, 2, 3, 4, 5, 4, 2, 3, 1, 4, 3 ] \\
\hline\phantom{\rule{.2pt}{10pt}}
	[ 1, 3, 1, 4, 2, 3, 4, 5, 4, 2, 3, 1, 4, 3 ] \\
\hline\phantom{\rule{.2pt}{10pt}}
	[ 1, 3, 4, 2, 3, 5, 4, 2, 3, 1, 6, 5, 4, 2, 3 ] \\
\hline\phantom{\rule{.2pt}{10pt}}
	[ 1, 3, 1, 4, 2, 3, 5, 4, 2, 3, 4, 6, 5, 4, 2, 3 ] \\
\hline\phantom{\rule{.2pt}{10pt}}
	[ 3, 1, 4, 2, 3, 5, 4, 2, 3, 1, 6, 5, 4, 2, 3, 4, 5 ] \\
\hline\phantom{\rule{.2pt}{10pt}}
	[ 3, 1, 4, 2, 3, 5, 4, 2, 3, 1, 6, 5, 4, 2, 3, 1, 4, 5 ] \\
\hline\phantom{\rule{.2pt}{10pt}}
	[ 3, 1, 4, 2, 3, 4, 5, 4, 2, 3, 4, 6, 5, 4, 2, 3, 1, 4, 3 ] \\
    \hline
\end{tabular}
\end{center}
}

%\end{document}

\vspace{5mm}
\paragraph{Table 4:  Conjugacy Classes Associated with $e_2$}
{
\footnotesize
\begin{center}
\begin{tabular}{|l|}
\hline\phantom{\rule{.2pt}{10pt}}
Conjugacy Classes  \\
	\hline
    [  ]	 \\
	\hline
	[ 4 ]	 \\
	\hline
	[ 3 ]	 \\
	\hline
	[ 3, 1 ] \\
	\hline
	[ 1, 4 ] \\
	\hline
	[ 4, 3 ] \\
	\hline
	[ 4, 3, 1 ] \\
	\hline
	[ 4, 2, 5, 4 ] \\
	\hline
	[ 4, 2, 3, 4 ] \\
	\hline
	[ 4, 3, 5, 4 ] \\
	\hline
	[ 4, 2, 3, 1, 4 ] \\
	\hline
	[ 4, 3, 1, 5, 4 ] \\
	\hline
	[ 4, 2, 5, 4, 3 ] \\
	\hline
	[ 4, 2, 3, 5, 4 ]\\
	\hline
	[ 1, 4, 2, 5, 4 ]\\
	\hline
	[ 4, 2, 5, 4, 3, 1 ]\\
	\hline
	[ 4, 2, 3, 1, 5, 4 ] \\
	\hline
	[ 4, 2, 3, 5, 4, 3 ] \\
	\hline
	[ 4, 2, 3, 1, 5, 4, 3 ] \\
	\hline
	[ 4, 2, 3, 5, 4, 3, 1 ] \\
	\hline
	[ 4, 3, 1, 5, 4, 2, 3, 4 ] \\
	\hline
	[ 4, 2, 3, 1, 5, 4, 3, 1 ] \\
	\hline
	[ 4, 2, 3, 1, 4, 3, 5, 4 ] \\
	\hline
	[ 4, 2, 5, 4, 3, 6, 5, 4 ] \\
	\hline
	[ 4, 3, 5, 4, 2, 6, 5, 4 ] \\
	\hline
	[ 4, 2, 3, 4, 5, 4, 2, 3, 4 ] \\
	\hline
	[ 4, 2, 3, 1, 5, 4, 2, 3, 4 ] \\
	\hline
	[ 4, 2, 3, 1, 4, 3, 5, 4, 3 ] \\
	\hline
	[ 4, 3, 1, 5, 4, 3, 6, 5, 4 ] \\
	\hline
	[ 4, 2, 5, 4, 3, 1, 6, 5, 4 ] \\
	\hline
\end{tabular}
\quad
\begin{tabular}{|l|}
\hline\phantom{\rule{.2pt}{10pt}}
Conjugacy Classes \\
	\hline
	[ 4, 3, 1, 5, 4, 2, 6, 5, 4 ] \\
	\hline
	[ 4, 2, 5, 4, 3, 1, 6, 5, 4, 3 ] \\
	\hline
	[ 4, 2, 3, 1, 5, 4, 3, 6, 5, 4 ] \\
	\hline
	[ 4, 2, 3, 4, 5, 4, 2, 3, 1, 4 ] \\
	\hline
	[ 4, 3, 5, 4, 2, 6, 5, 4, 3, 1 ]\\
	\hline
	[ 3, 4, 2, 3, 4, 5, 4, 2, 3, 4 ] \\
	\hline
	[ 4, 2, 3, 1, 4, 5, 4, 2, 3, 4, 6 ]\\
	\hline
	[ 3, 1, 4, 2, 3, 4, 5, 4, 2, 3, 4 ]\\
	\hline
	[ 4, 2, 3, 4, 5, 4, 2, 3, 6, 5, 4 ]\\
	\hline
	[ 3, 1, 4, 2, 3, 1, 4, 5, 4, 2, 3, 4 ] \\
	\hline
	[ 4, 2, 3, 1, 5, 4, 3, 1, 6, 5, 4, 3 ] \\
	\hline
	[ 4, 2, 3, 1, 4, 5, 4, 2, 3, 1, 4, 6 ] \\
	\hline
	[ 4, 2, 3, 4, 5, 4, 2, 3, 1, 6, 5, 4 ] \\
	\hline
	[ 4, 2, 3, 1, 4, 5, 4, 2, 3, 6, 5, 4 ] \\
	\hline
	[ 3, 4, 3, 1, 5, 4, 2, 3, 4, 6, 5, 4 ] \\
	\hline
	[ 3, 4, 2, 3, 1, 4, 5, 4, 2, 3, 4, 6 ] \\
	\hline
	[ 4, 2, 3, 4, 5, 4, 2, 3, 6, 5, 4, 3 ] \\
	\hline
	[ 4, 2, 3, 4, 5, 4, 2, 3, 1, 6, 5, 4, 3 ] \\
	\hline
	[ 4, 2, 3, 1, 4, 5, 4, 2, 3, 1, 6, 5, 4 ] \\
	\hline
	[ 1, 3, 4, 2, 3, 1, 5, 4, 2, 3, 6, 5, 4 ] \\
	\hline
	[ 3, 1, 4, 2, 3, 5, 4, 2, 3, 6, 5, 4, 2 ] \\
	\hline
	[ 4, 2, 3, 1, 4, 5, 4, 2, 3, 6, 5, 4, 3 ] \\
	\hline
	[ 1, 3, 1, 4, 3, 1, 5, 4, 2, 3, 6, 5, 4 ] \\
	\hline
	[ 1, 4, 2, 5, 4, 3, 1, 6, 5, 4, 2, 3, 4 ] \\
	\hline
	[ 1, 4, 2, 3, 1, 4, 3, 5, 4, 2, 6, 5, 4 ] \\
\hline
	[ 4, 2, 3, 4, 5, 4, 2, 3, 1, 4, 3, 6, 5, 4 ] \\
	\hline
	[ 4, 2, 3, 4, 5, 4, 3, 1, 6, 5, 4, 2, 3, 4 ] \\
	\hline
	[ 1, 4, 2, 3, 1, 4, 3, 5, 4, 2, 3, 6, 5, 4 ]\\
	\hline
	[ 3, 1, 4, 2, 3, 1, 5, 4, 2, 3, 6, 5, 4, 2 ] \\
	\hline
	[ 3, 1, 4, 2, 3, 5, 4, 2, 3, 1, 6, 5, 4, 2 ]\\
	\hline
\end{tabular}
\end{center}
}

\newpage
\paragraph{Table 4: Conjugacy Classes Associated with $e_2$ (continued)}
{
\footnotesize
\begin{center}
\begin{tabular}{|l|}
\hline\phantom{\rule{.2pt}{10pt}}
Conjugacy Classes  \\
	\hline
	[ 1, 3, 4, 2, 3, 1, 5, 4, 2, 3, 6, 5, 4, 3 ] \\
	\hline
	[ 4, 2, 3, 1, 4, 3, 5, 4, 2, 3, 1, 4, 3, 5, 4 ]\\
	\hline
	[ 4, 3, 1, 5, 4, 2, 3, 1, 4, 6, 5, 4, 2, 3, 4 ]\\
	\hline
	[ 4, 2, 3, 4, 5, 4, 2, 3, 1, 4, 3, 6, 5, 4, 3 ]\\
	\hline
	[ 3, 4, 2, 3, 1, 4, 3, 5, 4, 2, 3, 6, 5, 4, 2 ]\\
	\hline
	[ 4, 2, 3, 4, 5, 4, 2, 3, 1, 6, 5, 4, 2, 3, 4 ]\\
	\hline
	[ 3, 1, 4, 2, 3, 4, 5, 4, 2, 3, 1, 4, 3, 6, 5, 4 ]\\
\hline
\end{tabular}
\end{center}

}

\textheight 235mm
\topskip 0pt
\topmargin -12mm

\newpage
\paragraph{Table 5: Conjugacy Classes Associated with $e_3$}
{
\footnotesize
\begin{center}
\begin{tabular}{|l|}
\hline\phantom{\rule{.2pt}{10pt}}
Conjugacy Classes \\
	\hline
	[  ] \\
	\hline
	[ 5 ]\\
	\hline
	[ 4 ]\\
	\hline
	[ 2 ]\\
	\hline
	[ 4, 2 ]\\
	\hline
	[ 2, 5 ]\\
	\hline
	[ 5, 4 ]\\
	\hline
	[ 2, 3 ]\\
	\hline
	[ 3, 1 ]\\
	\hline
	[ 1, 4 ]\\
	\hline
	[ 3, 5 ]\\
	\hline
	[ 5, 4, 2 ]\\
	\hline
	[ 4, 2, 5 ]\\
	\hline
	[ 2, 4, 5 ]\\
	\hline
	[ 1, 4, 2 ]\\
	\hline
	[ 4, 2, 3 ]\\
	\hline
	[ 2, 3, 1 ]\\
	\hline
	[ 3, 1, 5 ]\\
	\hline
	[ 4, 3, 1 ]\\
	\hline
	[ 5, 4, 3 ]\\
	\hline
	[ 1, 5, 4 ]\\
	\hline
	[ 2, 3, 5 ]\\
	\hline
	[ 4, 5, 4, 2 ]\\
	\hline
	[ 5, 4, 6, 5 ]\\
	\hline
	[ 2, 4, 2, 5 ]\\
	\hline
	[ 4, 2, 3, 1 ]\\
	\hline
	[ 1, 5, 4, 2 ]\\
	\hline
	[ 4, 2, 3, 5 ]\\
	\hline
	[ 5, 4, 3, 1 ]\\
	\hline
	[ 2, 3, 1, 5 ]\\
	\hline
	[ 1, 4, 2, 5 ]\\
	\hline
	[ 1, 2, 4, 5 ]\\
	\hline
	[ 5, 4, 2, 3 ]\\
	\hline
	[ 2, 4, 3, 5 ]\\
	\hline
	[ 2, 4, 5, 4, 2 ]\\
	\hline
	[ 5, 4, 2, 6, 5 ]\\
	\hline
	[ 2, 5, 4, 6, 5 ]\\
	\hline
	[ 4, 5, 4, 2, 3 ]\\
	\hline
	[ 2, 4, 5, 4, 3 ]\\
	\hline
	[ 1, 4, 5, 4, 2 ]\\
	\hline
	[ 5, 4, 2, 3, 1 ]\\
	\hline
	[ 4, 2, 3, 1, 5 ]\\
    \hline
	[ 2, 4, 3, 1, 5 ]\\
	\hline
	[ 1, 2, 4, 2, 5 ]\\
	\hline
	[ 5, 4, 3, 6, 5 ]\\
	\hline
	[ 2, 5, 4, 2, 3 ]\\
	\hline
	[ 1, 5, 4, 6, 5 ]\\
	\hline
	[ 2, 4, 2, 3, 5 ]\\
	\hline
	[ 5, 4, 2, 3, 4, 5 ]\\
	\hline
	[ 2, 5, 4, 2, 6, 5 ]\\
	\hline
	[ 1, 2, 4, 5, 4, 2 ]\\
	\hline
	[ 2, 4, 3, 5, 4, 2 ]\\
	\hline
	[ 2, 4, 2, 5, 4, 2 ]\\
	\hline
	[ 2, 5, 4, 2, 3, 1 ]\\
	\hline
\end{tabular}   \quad
\begin{tabular}{|l|}
\hline\phantom{\rule{.2pt}{10pt}}
Conjugacy Classes \\
	\hline
	[ 4, 5, 4, 2, 3, 1 ] \\
	\hline
	[ 4, 2, 5, 4, 2, 3 ]\\
	\hline
	[ 2, 4, 2, 5, 4, 3 ]\\
	\hline
	[ 2, 4, 5, 4, 3, 1 ]\\
	\hline
	[ 5, 4, 3, 1, 6, 5 ]\\
	\hline
	[ 2, 4, 2, 3, 1, 5 ]\\
	\hline
	[ 1, 2, 5, 4, 6, 5 ]\\
\hline
\end{tabular}
\quad
\begin{tabular}{|l|}
\hline\phantom{\rule{.2pt}{10pt}}
Conjugacy Classes \\
	\hline
	[ 2, 5, 4, 2, 6, 5, 4, 3 ]\\
	\hline
	[ 4, 5, 4, 2, 3, 1, 4, 5 ]\\
	\hline
	[ 5, 4, 2, 3, 4, 6, 5, 4 ]\\
	\hline
	[ 2, 4, 2, 3, 1, 5, 4, 2 ]\\
	\hline
	[ 1, 4, 2, 5, 4, 2, 6, 5 ]\\
	\hline
	[ 2, 4, 3, 5, 4, 2, 3, 4, 5 ]\\
	\hline
	[ 3, 1, 5, 4, 2, 3, 4, 6, 5 ]\\
\hline
\end{tabular}
\end{center}
}

\newpage
\paragraph{Table 5: Conjugacy Classes Associated with $e_3$ (continued)}

{
\footnotesize
\begin{center}
\begin{tabular}{|l|}
\hline\phantom{\rule{.2pt}{10pt}}
Conjugacy Classes \\
	\hline
	[ 2, 4, 2, 3, 5, 4, 3, 1, 6, 5 ]\\
	\hline
	[ 2, 4, 3, 1, 5, 4, 3, 6, 5, 4 ]\\
	\hline
	[ 2, 3, 1, 4, 3, 5, 4, 2, 6, 5 ]\\
	\hline
	[ 2, 4, 3, 5, 4, 3, 6, 5, 4, 2 ]\\
	\hline
	[ 2, 4, 2, 3, 1, 5, 4, 2, 6, 5 ]\\
	\hline
	[ 2, 4, 2, 3, 5, 4, 2, 3, 6, 5 ]	\\
	\hline
	[ 2, 4, 2, 3, 5, 4, 2, 3, 4, 5 ]	\\
\hline
\end{tabular}\quad
\begin{tabular}{|l|}
\hline\phantom{\rule{.2pt}{10pt}}
Conjugacy Classes \\
	\hline
	[ 4, 2, 3, 1, 4, 3, 5, 4, 2, 3, 4, 5 ]\\
	\hline
	[ 4, 2, 5, 4, 2, 3, 4, 6, 5, 4, 2, 3 ]	\\
	\hline
	[ 4, 2, 3, 1, 5, 4, 3, 1, 6, 5, 4, 2 ]\\
	\hline
	[ 5, 4, 2, 3, 1, 6, 5, 4, 2, 3, 4, 5 ]\\
	\hline
	[ 2, 4, 5, 4, 2, 3, 1, 4, 6, 5, 4, 2 ]\\
	\hline
	[ 4, 2, 3, 1, 5, 4, 2, 3, 1, 4, 6, 5 ]\\
	\hline
	[ 4, 3, 1, 5, 4, 2, 3, 1, 4, 3, 6, 5 ]\\
\hline
\end{tabular}
\end{center}
}

\newpage
\paragraph{Table 5: Conjugacy Classes Associated with $e_3$ (continued)}
{
\footnotesize
\begin{center}
\begin{tabular}{|l|}
\hline\phantom{\rule{.2pt}{10pt}}
Conjugacy Classes \\
	\hline
	[ 2, 4, 5, 4, 3, 1, 6, 5, 4, 2, 3, 4, 5 ]	\\
	\hline
	[ 3, 1, 5, 4, 2, 3, 4, 6, 5, 4, 2, 3, 4 ]	\\
	\hline
	[ 3, 1, 4, 2, 3, 4, 5, 4, 2, 3, 4, 6, 5 ]	\\
	\hline
	[ 2, 3, 1, 4, 2, 3, 4, 5, 4, 2, 3, 4, 5 ]	\\
	\hline
	[ 3, 1, 5, 4, 2, 3, 1, 4, 3, 6, 5, 4, 2 ]	\\
	\hline
	[ 2, 4, 5, 4, 2, 3, 1, 4, 3, 6, 5, 4, 2 ]	\\
	\hline
	[ 2, 4, 2, 3, 5, 4, 2, 3, 1, 4, 3, 6, 5 ]	\\
\hline
\end{tabular}\quad
\begin{tabular}{|l|}
\hline\phantom{\rule{.2pt}{10pt}}
Conjugacy Classes \\
	\hline
	[ 2, 4, 3, 1, 5, 4, 2, 3, 4, 6, 5, 4, 2, 3 ]\\
	\hline
	[ 2, 4, 5, 4, 2, 3, 1, 4, 3, 5, 4, 2, 6, 5 ]\\
	\hline
	[ 2, 4, 2, 5, 4, 2, 3, 1, 4, 3, 6, 5, 4, 2 ]\\
	\hline
	[ 4, 3, 1, 5, 4, 2, 3, 1, 4, 3, 6, 5, 4, 2 ]\\
	\hline
	[ 4, 5, 4, 2, 3, 1, 4, 3, 6, 5, 4, 2, 3, 1 ]\\
	\hline
	[ 2, 4, 2, 3, 1, 5, 4, 3, 1, 6, 5, 4, 2, 3 ]\\
	\hline
	[ 5, 4, 2, 3, 1, 4, 3, 5, 4, 6, 5, 4, 2, 3 ]\\
\hline
\end{tabular}
\end{center}
}

\newpage
\paragraph{Table 5: Conjugacy Classes Associated with $e_3$ (continued)}
{
\footnotesize
\begin{center}
\begin{tabular}{|l|}
\hline\phantom{\rule{.2pt}{10pt}}
Conjugacy Classes \\
    \hline
	[ 2, 4, 2, 3, 1, 5, 4, 2, 3, 1, 4, 3, 5, 4, 2 ]\\
	\hline
	[ 2, 4, 3, 1, 5, 4, 2, 3, 1, 4, 3, 6, 5, 4, 2 ]\\
	\hline
	[ 2, 4, 3, 1, 5, 4, 3, 1, 6, 5, 4, 2, 3, 4, 5 ]\\
	\hline
	[ 4, 2, 5, 4, 2, 3, 1, 4, 3, 6, 5, 4, 2, 3, 1 ]\\
	\hline
	[ 2, 4, 2, 5, 4, 2, 3, 1, 4, 3, 6, 5, 4, 2, 3 ]\\
	\hline
	[ 5, 4, 2, 3, 1, 4, 3, 5, 4, 2, 6, 5, 4, 2, 3 ]\\
	\hline
	[ 4, 2, 3, 1, 4, 3, 5, 4, 3, 1, 6, 5, 4, 2, 3 ]\\
\hline
\end{tabular}
\quad
\begin{tabular}{|l|}
\hline\phantom{\rule{.2pt}{10pt}}
Conjugacy Classes \\
	\hline
	[ 2, 4, 2, 5, 4, 2, 3, 1, 4, 3, 6, 5, 4, 2, 3, 1 ]\\
	\hline
	[ 5, 4, 2, 3, 1, 4, 3, 5, 4, 2, 6, 5, 4, 2, 3, 4 ]\\
	\hline
	[ 2, 4, 2, 3, 1, 4, 3, 5, 4, 3, 1, 6, 5, 4, 2, 3 ]\\
	\hline
	[ 2, 4, 2, 3, 5, 4, 2, 3, 1, 4, 3, 6, 5, 4, 2, 3 ]\\
	\hline
	[ 4, 2, 5, 4, 2, 3, 1, 4, 5, 6, 5, 4, 2, 3, 4, 5 ]\\
	\hline
	[ 5, 4, 2, 3, 1, 4, 3, 5, 4, 6, 5, 4, 2, 3, 4, 5 ]\\
	\hline
	[ 2, 5, 4, 2, 3, 1, 4, 5, 6, 5, 4, 2, 3, 1, 4, 5 ]\\
\hline
\end{tabular}
\end{center}
}

\newpage
\paragraph{Table 5: Conjugacy Classes Associated with $e_3$ (continued)}
{
\footnotesize

\begin{center}
\begin{tabular}{|l|}
\hline\phantom{\rule{.2pt}{10pt}}
Conjugacy Classes \\
	\hline
	[ 2, 5, 4, 2, 3, 1, 4, 3, 6, 5, 4, 2, 3, 1, 4, 3, 5, 4 ]\\
	\hline
	[ 2, 4, 2, 3, 1, 5, 4, 2, 3, 4, 6, 5, 4, 2, 3, 1, 4, 3 ]\\
	\hline
	[ 2, 4, 3, 5, 4, 2, 3, 1, 4, 3, 5, 4, 2, 6, 5, 4, 2, 3 ]\\
	\hline
	[ 2, 4, 2, 3, 5, 4, 2, 3, 1, 4, 3, 5, 4, 6, 5, 4, 2, 3 ]\\
	\hline
	[ 2, 4, 5, 4, 2, 3, 1, 4, 3, 5, 6, 5, 4, 2, 3, 1, 4, 5 ]\\
	\hline
	[ 5, 4, 2, 3, 1, 4, 3, 5, 6, 5, 4, 2, 3, 1, 4, 3, 5, 4 ]\\
	\hline
	[ 4, 3, 1, 5, 4, 2, 3, 1, 4, 3, 5, 6, 5, 4, 2, 3, 4, 5 ]\\
\hline
\end{tabular}
\end{center}
}

\newpage
\paragraph{Table 5: Conjugacy Classes Associated with $e_3$ (continued)}

{
\footnotesize
\begin{center}
\begin{tabular}{|l|}
\hline\phantom{\rule{.2pt}{10pt}}
Conjugacy Classes \\
	\hline
	[ 2, 4, 2, 3, 1, 5, 4, 2, 3, 1, 4, 5, 6, 5, 4, 2, 3, 4, 5 ] \\
	\hline
	[ 2, 4, 5, 4, 2, 3, 1, 4, 3, 6, 5, 4, 2, 3, 1, 4, 3, 5, 4 ] \\
	\hline
	[ 2, 5, 4, 2, 3, 1, 4, 3, 5, 4, 6, 5, 4, 2, 3, 1, 4, 3, 5 ] \\
	\hline
	[ 2, 4, 3, 1, 5, 4, 2, 3, 1, 4, 3, 6, 5, 4, 2, 3, 1, 4, 5 ] \\
	\hline
	[ 2, 3, 1, 4, 5, 4, 2, 3, 1, 4, 3, 5, 6, 5, 4, 2, 3, 4, 5 ] \\
	\hline
	[ 2, 4, 5, 4, 2, 3, 1, 4, 3, 5, 4, 2, 6, 5, 4, 2, 3, 4, 5 ] \\
	\hline
	[ 3, 4, 2, 3, 1, 4, 3, 5, 4, 2, 3, 1, 4, 3, 5, 4, 2, 6, 5 ] \\
\hline
\end{tabular}
\end{center}
}

\newpage
\paragraph{Table 5: Conjugacy Classes Associated with $e_3$ (continued)}

{
\footnotesize
\begin{center}
\begin{tabular}{|l|}
\hline\phantom{\rule{.2pt}{10pt}}
Conjugacy Classes \\
	\hline	
	[ 2, 4, 3, 1, 5, 4, 2, 3, 4, 5, 6, 5, 4, 2, 3, 1, 4, 3, 5, 4, 2 ]
\\
	\hline
	[ 2, 4, 2, 5, 4, 2, 3, 1, 4, 3, 5, 4, 6, 5, 4, 2, 3, 1, 4, 3, 5 ]
\\
	\hline
	[ 3, 1, 4, 3, 5, 4, 2, 3, 1, 4, 3, 5, 4, 6, 5, 4, 2, 3, 1, 4, 5 ]
\\
    \hline
    [ 3, 1, 4, 2, 3, 4, 5, 4, 2, 3, 1, 4, 3, 5, 6, 5, 4, 2, 3, 4, 5 ]
\\
	\hline
	[ 2, 3, 1, 4, 3, 5, 4, 2, 3, 1, 4, 3, 5, 4, 6, 5, 4, 2, 3, 4, 5 ]
\\
	\hline
	[ 2, 3, 4, 3, 1, 5, 4, 2, 3, 4, 5, 6, 5, 4, 2, 3, 1, 4, 3, 5, 4 ]
\\
	\hline
	[ 2, 4, 3, 1, 5, 4, 2, 3, 1, 4, 3, 5, 4, 2, 6, 5, 4, 2, 3, 4, 5 ]
\\
\hline
\end{tabular}
\end{center}
}

\newpage
\paragraph{Table 5: Conjugacy Classes Associated with $e_3$ (continued)}

{
\footnotesize
\begin{center}
\begin{tabular}{|l|}
\hline\phantom{\rule{.2pt}{10pt}}
Conjugacy Classes \\
	\hline
	[ 3, 1, 4, 2, 3, 5, 4, 2, 3, 1, 4, 5, 6, 5, 4, 2, 3, 1, 4, 3, 5, 4, 2
] \\
	\hline
	[ 3, 5, 4, 2, 3, 1, 4, 3, 5, 4, 2, 6, 5, 4, 2, 3, 1, 4, 3, 5, 4, 6, 5
] \\
	\hline
	[ 3, 1, 4, 3, 1, 5, 4, 2, 3, 1, 4, 5, 6, 5, 4, 2, 3, 1, 4, 3, 5, 4, 2
] \\
	\hline
	[ 2, 3, 1, 4, 3, 1, 5, 4, 2, 3, 4, 5, 6, 5, 4, 2, 3, 1, 4, 3, 5, 4, 2
] \\
	\hline
	[ 2, 3, 1, 4, 3, 1, 5, 4, 2, 3, 1, 4, 5, 6, 5, 4, 2, 3, 1, 4, 3, 5, 4
] \\
	\hline
	[ 2, 4, 5, 4, 2, 3, 1, 4, 3, 5, 4, 6, 5, 4, 2, 3, 1, 4, 3, 5, 4, 6, 5
] \\
	\hline
	[ 1, 2, 4, 2, 3, 1, 5, 4, 2, 3, 4, 5, 6, 5, 4, 2, 3, 1, 4, 3, 5, 4, 2
] \\
\hline
\end{tabular}
\end{center}

}

\textheight 240mm
\topskip 0pt
\topmargin -12mm

\newpage
\paragraph{Table 6: Conjugacy Classes Associated with $e_4$}

{
\footnotesize

\begin{center}
\begin{tabular}{|l|}
\hline\phantom{\rule{.2pt}{10pt}}
Conjugacy Classes \\
	\hline
    [  ]	 \\
	\hline
	[ 5 ]	 \\
	\hline
	[ 4 ]	 \\
	\hline
	[ 4, 2 ]	 \\
	\hline
	[ 2, 5 ]	 \\
	\hline
	[ 5, 4 ]	 \\
	\hline
	[ 3, 5 ]	 \\
	\hline
	[ 2, 3 ]	 \\
	\hline
	[ 5, 4, 2 ]	 \\
	\hline
	[ 3, 1, 5 ]	 \\
	\hline
	[ 2, 3, 1 ]	 \\
	\hline
	[ 1, 5, 4 ]	\\
	\hline
	[ 4, 3, 1 ]	 \\
	\hline
	[ 5, 4, 3 ]	\\
	\hline
	[ 2, 3, 5 ]	 \\
	\hline
	[ 5, 4, 6, 5 ]\\
	\hline
	[ 4, 2, 3, 1 ]	 \\
	\hline
	[ 5, 4, 3, 1 ]	 \\
	\hline
	[ 1, 5, 4, 2 ]	\\
	\hline
	[ 5, 4, 2, 3 ]	 \\
	\hline
	[ 2, 3, 1, 5 ]	 \\
	\hline
	[ 5, 4, 2, 6, 5 ]\\
	\hline
	[ 4, 5, 4, 2, 3 ]	\\
	\hline
	[ 5, 4, 2, 3, 1 ]	 \\
	\hline
	[ 5, 4, 3, 6, 5 ]	 \\
	\hline
	[ 1, 5, 4, 6, 5 ]	 \\
	\hline
	[ 5, 4, 2, 3, 4, 5 ]	\\
	\hline
	[ 4, 5, 4, 2, 3, 1 ]	 \\
	\hline
	[ 4, 2, 5, 4, 2, 3 ]	 \\
	\hline
	[ 5, 4, 3, 1, 6, 5 ]\\
	\hline
	[ 1, 5, 4, 2, 6, 5 ]	 \\
	\hline
	[ 5, 4, 2, 3, 6, 5 ]	\\
	\hline
	[ 4, 5, 4, 2, 3, 4, 5 ]	 \\
	\hline
	[ 5, 4, 2, 3, 1, 4, 5 ]	 \\
	\hline
	[ 5, 4, 2, 3, 4, 6, 5 ]	 \\
	\hline
	[ 3, 1, 5, 4, 2, 3, 4 ]	 \\
	\hline
	[ 2, 5, 4, 2, 3, 4, 5 ]	 \\
	\hline
	[ 5, 4, 2, 3, 1, 6, 5 ]	 \\
	\hline
	[ 4, 2, 5, 4, 3, 6, 5 ]	\\
	\hline
	[ 5, 4, 2, 6, 5, 4, 3 ]	\\
	\hline
	[ 3, 1, 5, 4, 2, 3, 4, 5 ]	 \\
	\hline
	[ 5, 4, 3, 1, 6, 5, 4, 3 ]\\
	\hline
	[ 5, 4, 2, 3, 1, 4, 6, 5 ]	 \\
	\hline
	[ 2, 5, 4, 2, 3, 4, 6, 5 ]	\\
	\hline
	[ 5, 4, 2, 6, 5, 4, 3, 1 ]	\\
	\hline
	[ 5, 4, 2, 3, 4, 6, 5, 4 ]	 \\
	\hline
	[ 4, 2, 5, 4, 3, 1, 6, 5 ]	 \\
	\hline
	[ 4, 5, 4, 2, 3, 1, 4, 5 ]	 \\
	\hline
	[ 4, 2, 5, 4, 2, 3, 4, 5 ]	 \\
	\hline
	[ 3, 1, 5, 4, 2, 3, 4, 6, 5 ]	 \\
	\hline
	[ 5, 4, 2, 3, 1, 6, 5, 4, 3 ]	 \\
	\hline
	[ 4, 3, 1, 5, 4, 2, 3, 4, 5 ]	 \\
	\hline
	[ 2, 5, 4, 3, 1, 6, 5, 4, 3 ]	 \\
	\hline
	[ 5, 4, 2, 3, 1, 4, 6, 5, 4 ]	 \\
	\hline
	[ 5, 4, 2, 3, 4, 6, 5, 4, 3 ]	 \\
	\hline
	[ 4, 2, 5, 4, 2, 3, 4, 6, 5 ]	 \\
	\hline
	[ 4, 2, 5, 4, 2, 6, 5, 4, 3 ]	 \\
	\hline
\end{tabular}
\quad
\begin{tabular}{|l|}
\hline\phantom{\rule{.2pt}{10pt}}
Conjugacy Classes \\
	\hline
	[ 4, 3, 1, 5, 4, 2, 3, 1, 4, 5 ]	 \\
	\hline
	[ 2, 5, 4, 2, 3, 1, 6, 5, 4, 3 ]	\\
	\hline
	[ 5, 4, 2, 3, 1, 4, 3, 6, 5, 4 ]	 \\
	\hline
	[ 4, 3, 1, 5, 4, 2, 3, 4, 6, 5 ]	 \\
	\hline
	[ 4, 2, 5, 4, 2, 6, 5, 4, 3, 1 ]	\\
	\hline
	[ 5, 4, 2, 3, 1, 4, 6, 5, 4, 3 ]	 \\
	\hline
	[ 4, 2, 5, 4, 2, 3, 1, 4, 6, 5 ]	 \\
\hline
\end{tabular}
\end{center}
}

\newpage
\paragraph{Table 6: Conjugacy Classes Associated with $e_4$ (continued)}

{
\footnotesize

\begin{center}
\begin{tabular}{|l|}
\hline\phantom{\rule{.2pt}{10pt}}
Conjugacy Classes \\
	\hline
	[ 2, 3, 5, 4, 2, 3, 4, 5, 6, 5, 4, 2, 3, 4, 5 ]	 \\
	\hline
	[ 4, 2, 5, 4, 2, 3, 1, 6, 5, 4, 2, 3, 1, 4, 5 ]	 \\
	\hline
	[ 5, 4, 2, 3, 1, 4, 3, 5, 4, 2, 6, 5, 4, 2, 3 ]	 \\
	\hline
	[ 4, 5, 4, 2, 3, 4, 5, 6, 5, 4, 2, 3, 1, 4, 5 ]	 \\
	\hline
	[ 3, 1, 4, 3, 5, 4, 2, 3, 1, 4, 3, 5, 4, 6, 5 ]	 \\
	\hline
	[ 3, 1, 5, 4, 2, 3, 1, 4, 3, 6, 5, 4, 2, 3, 4 ]	 \\
	\hline
	[ 4, 5, 4, 2, 3, 1, 4, 3, 5, 4, 2, 6, 5, 4, 3 ]	 \\
\hline
\end{tabular}
\end{center}

}

\textheight 235mm
\topskip 0pt
\topmargin -12mm

\newpage
\paragraph{Table 7: Conjugacy Classes Associated with $e_5$}

{
\footnotesize

\begin{center}
\begin{tabular}[t]{|l|}
\hline\phantom{\rule{.2pt}{10pt}}
Conjugacy Classes\\
	\hline
	[  ]	\\
	\hline
	[ 6 ]	\\
	\hline
	[ 5 ]	\\
	\hline
	[ 6, 5 ]	\\
	\hline
	[ 3, 5 ]	\\
	\hline
	[ 4, 2 ]	\\
	\hline
	[ 2, 5 ]	\\
	\hline
	[ 4, 6 ]	\\
	\hline
	[ 5, 4, 2 ]	\\
	\hline
	[ 6, 5, 4 ]	\\
	\hline
	[ 4, 2, 6 ]	\\
	\hline
	[ 3, 1, 5 ]	\\
	\hline
	[ 2, 6, 5 ]	\\
	\hline
	[ 2, 3, 5 ]	\\
	\hline
	[ 3, 6, 5 ]	\\
	\hline
	[ 4, 3, 1 ]	\\
	\hline
	[ 2, 3, 6 ]	\\
	\hline
	[ 6, 5, 4, 2 ]	\\
	\hline
	[ 1, 5, 4, 2 ]	\\
	\hline
	[ 3, 1, 6, 5 ]	\\
	\hline
	[ 4, 2, 3, 1 ]	\\
	\hline
	[ 1, 6, 5, 4 ]	\\
	\hline
	[ 2, 3, 1, 6 ]	\\
	\hline
	[ 6, 5, 4, 3 ]	\\
	\hline
	[ 4, 3, 1, 6 ]	\\
	\hline
	[ 5, 4, 2, 3 ]	\\
	\hline
	[ 2, 3, 6, 5 ]	\\
	\hline
	[ 2, 3, 1, 5 ]	\\
	\hline
	[ 4, 2, 3, 1, 6 ]\\
	\hline
	[ 6, 5, 4, 3, 1 ]\\
	\hline
	[ 1, 6, 5, 4, 2 ]\\
	\hline
	[ 5, 4, 2, 3, 1 ]\\
	\hline
	[ 6, 5, 4, 2, 3 ]\\
	\hline
	[ 2, 3, 1, 6, 5 ]\\
	\hline
	[ 5, 6, 5, 4, 2, 3 ]\\
	\hline
	[ 6, 5, 4, 2, 3, 1 ]\\
	\hline
	[ 4, 2, 5, 4, 2, 3 ]\\
	\hline
	[ 4, 5, 6, 5, 4, 2, 3 ]\\
	\hline
	[ 3, 1, 5, 4, 2, 3, 4 ]\\
	\hline
	[ 5, 6, 5, 4, 2, 3, 1 ]\\
	\hline
	[ 5, 4, 2, 6, 5, 4, 3 ]\\
	\hline
	[ 6, 5, 4, 2, 3, 4, 5, 6 ]\\
	\hline
	[ 4, 5, 6, 5, 4, 2, 3, 1 ]\\
	\hline

\end{tabular}
\quad
\begin{tabular}[t]{|l|}
\hline\phantom{\rule{.2pt}{10pt}}
Conjugacy Classes \\
	\hline	
	[ 4, 5, 4, 6, 5, 4, 2, 3 ]\\
	\hline
	[ 3, 1, 6, 5, 4, 2, 3, 4 ]\\
	\hline
    [ 4, 2, 5, 6, 5, 4, 2, 3 ]\\
	\hline
	[ 2, 6, 5, 4, 2, 3, 4, 5, 6 ]\\
	\hline
	[ 5, 6, 5, 4, 2, 3, 4, 5, 6 ]\\
	\hline
	[ 6, 5, 4, 2, 3, 1, 4, 5, 6 ]\\
	\hline
	[ 3, 1, 5, 6, 5, 4, 2, 3, 4 ]\\
\hline
\end{tabular}
\end{center}

}

\newpage
\paragraph{Table 8: Conjugacy Classes Associated with $e_6$}

{
\footnotesize

\begin{center}
\begin{tabular}{|l|}
\hline\phantom{\rule{.2pt}{10pt}}
Conjugacy Classes \\
	\hline
	[  ]	\\
	\hline
	[ 6 ]	\\
	\hline
	[ 5 ]	\\
	\hline
	[ 2 ]	\\
	\hline
	[ 6, 5 ]	 \\
	\hline
	[ 2, 3 ]	 \\
	\hline
	[ 4, 2 ]	 \\
	\hline
	[ 2, 5 ]	 \\
	\hline
	[ 2, 6 ]	 \\
	\hline
	[ 5, 4 ]	 \\
	\hline
	[ 4, 6 ]	\\
	\hline
	[ 3, 5 ]	 \\
	\hline
	[ 5, 4, 2 ]	 \\
	\hline
	[ 4, 2, 6 ]	 \\
	\hline
	[ 2, 6, 5 ]	 \\
	\hline
	[ 6, 5, 4 ]	 \\
	\hline
	[ 2, 3, 1 ]	 \\
	\hline
	[ 4, 2, 3 ]	 \\
	\hline
	[ 2, 3, 6 ]	 \\
	\hline
	[ 1, 4, 2 ]	 \\
	\hline
	[ 3, 1, 6 ]	 \\
	\hline
	[ 3, 1, 5 ]	 \\
	\hline
	[ 4, 3, 1 ]	 \\
	\hline
	[ 1, 4, 6 ]	 \\
	\hline
	[ 3, 6, 5 ]	 \\
	\hline
	[ 2, 3, 5 ]	 \\
	\hline
	[ 6, 5, 4, 2 ]	 \\
	\hline
	[ 5, 4, 2, 6 ]	 \\
	\hline
	[ 4, 2, 5, 6 ]	 \\
	\hline
	[ 2, 4, 5, 6 ]	 \\
	\hline
	[ 4, 2, 3, 1 ]	 \\
	\hline
	[ 1, 4, 2, 6 ]	 \\
	\hline
	[ 1, 5, 4, 2 ]	 \\
	\hline
	[ 4, 2, 3, 6 ]	 \\
	\hline
	[ 2, 3, 1, 6 ]	 \\
	\hline
	[ 3, 1, 6, 5 ]	 \\
	\hline
	[ 5, 4, 2, 3 ]	 \\
	\hline
	[ 5, 4, 3, 1 ]	 \\
	\hline
	[ 6, 5, 4, 3 ]	\\
	\hline
	[ 4, 3, 1, 6 ]	 \\
	\hline
	[ 1, 6, 5, 4 ]	\\
	\hline
	[ 2, 3, 1, 5 ]	\\
	\hline
	[ 2, 3, 6, 5 ]	\\
	\hline
	[ 5, 6, 5, 4, 2 ]	 \\
	\hline
	[ 4, 5, 4, 2, 6 ]\\
	\hline
    [ 2, 4, 2, 5, 6 ]	 \\
	\hline
	[ 4, 2, 3, 1, 6 ]	 \\
	\hline
	[ 1, 6, 5, 4, 2 ]	 \\
	\hline
	[ 5, 4, 2, 3, 1 ]	 \\
	\hline
	[ 1, 5, 4, 2, 6 ]	 \\
	\hline
	[ 5, 4, 2, 3, 6 ]	 \\
	\hline
	[ 4, 2, 3, 5, 6 ]	 \\
	\hline
	[ 6, 5, 4, 3, 1 ]	 \\
	\hline
	[ 2, 4, 3, 5, 6 ]	\\
	\hline
	[ 1, 2, 4, 5, 6 ]	\\
	\hline
\end{tabular}   \quad
\begin{tabular}{|l|}
\hline\phantom{\rule{.2pt}{10pt}}
Conjugacy Classes \\
	\hline	
	[ 6, 5, 4, 2, 3 ]	\\
	\hline
	[ 2, 3, 1, 6, 5 ]	\\
	\hline
	[ 2, 5, 4, 2, 3 ]	 \\
	\hline
	[ 1, 4, 2, 5, 6 ]	 \\
	\hline
	[ 4, 5, 6, 5, 4, 2 ]	 \\
	\hline
	[ 2, 4, 5, 4, 2, 6 ]	 \\
	\hline
	[ 2, 4, 3, 5, 4, 2 ]	 \\
\hline
\end{tabular}
\quad
\begin{tabular}{|l|}
\hline\phantom{\rule{.2pt}{10pt}}
Conjugacy Classes \\
	\hline	
	[ 3, 1, 6, 5, 4, 2, 3, 4 ]	\\
	\hline
	[ 3, 1, 5, 4, 2, 3, 4, 6 ]	 \\
	\hline
	[ 3, 1, 4, 2, 3, 4, 5, 6 ]	\\
	\hline
	[ 2, 5, 6, 5, 4, 2, 3, 1 ]	 \\
	\hline
	[ 6, 5, 4, 2, 3, 1, 4, 3 ]	\\
	\hline
	[ 4, 2, 5, 6, 5, 4, 2, 3 ]	\\
	\hline
	[ 2, 4, 2, 5, 6, 5, 4, 3 ]	\\
\hline
\end{tabular}
\end{center}

}

\newpage
\paragraph{Table 8: Conjugacy Classes Associated with $e_6$ (continued)}

{
\footnotesize

\begin{center}
\begin{tabular}{|l|}
\hline\phantom{\rule{.2pt}{10pt}}
Conjugacy Classes \\
	\hline	
	[ 2, 4, 3, 5, 6, 5, 4, 2, 3 ]	 \\
	\hline
	[ 1, 2, 4, 2, 5, 6, 5, 4, 2 ]	 \\
	\hline
	[ 2, 4, 3, 1, 5, 4, 2, 3, 4, 5 ]	 \\
	\hline
	[ 4, 2, 5, 4, 3, 6, 5, 4, 2, 3 ]	 \\
	\hline
	[ 2, 4, 3, 5, 4, 2, 3, 4, 5, 6 ]	 \\
	\hline
	[ 6, 5, 4, 2, 3, 1, 4, 3, 5, 6 ]	 \\
	\hline
	[ 3, 1, 6, 5, 4, 2, 3, 4, 5, 6 ]	 \\
\hline
\end{tabular}
\quad
\begin{tabular}{|l|}
\hline\phantom{\rule{.2pt}{10pt}}
Conjugacy Classes \\
	\hline
	[ 2, 4, 5, 4, 2, 3, 6, 5, 4, 2, 3 ]	 \\
	\hline
	[ 5, 6, 5, 4, 2, 3, 1, 4, 3, 5, 6 ]	 \\
	\hline
	[ 2, 3, 1, 4, 2, 3, 4, 6, 5, 4, 2 ]	 \\
	\hline
	[ 2, 4, 2, 3, 1, 5, 4, 2, 3, 4, 5 ]	 \\
	\hline
	[ 2, 3, 1, 4, 2, 3, 4, 5, 4, 2, 6 ]	 \\
	\hline
	[ 2, 4, 3, 1, 6, 5, 4, 2, 3, 1, 4 ]	 \\
	\hline
	[ 2, 4, 3, 1, 5, 4, 2, 3, 1, 4, 6 ]	 \\
\hline
\end{tabular}
\end{center}

}

\newpage
\paragraph{Table 8: Conjugacy Classes Associated with $e_6$ (continued)}

{
\footnotesize

\begin{center}
\begin{tabular}{|l|}
\hline\phantom{\rule{.2pt}{10pt}}
Conjugacy Classes \\
	\hline
	[ 3, 1, 4, 2, 3, 1, 4, 5, 6, 5, 4, 2 ]	 \\
	\hline
	[ 2, 3, 1, 4, 3, 5, 4, 2, 3, 4, 5, 6 ]	 \\
	\hline
	[ 2, 4, 3, 1, 5, 4, 3, 1, 6, 5, 4, 2 ]	 \\
	\hline
	[ 4, 2, 5, 4, 2, 3, 4, 5, 6, 5, 4, 2 ]	 \\
	\hline
	[ 2, 4, 2, 5, 4, 2, 3, 6, 5, 4, 2, 3 ]	 \\
	\hline
	[ 4, 2, 3, 1, 4, 3, 5, 4, 2, 3, 6, 5 ]	 \\
	\hline
	[ 2, 4, 2, 3, 6, 5, 4, 2, 3, 4, 5, 6 ]	 \\
\hline
\end{tabular}
\quad
\begin{tabular}{|l|}
\hline\phantom{\rule{.2pt}{10pt}}
Conjugacy Classes \\
	\hline	
	[ 3, 1, 5, 4, 2, 3, 1, 4, 3, 5, 6, 5, 4, 2 ]	 \\
	\hline
	[ 4, 2, 5, 4, 3, 1, 6, 5, 4, 2, 3, 4, 5, 6 ]	 \\
	\hline
	[ 2, 4, 2, 3, 1, 5, 4, 3, 6, 5, 4, 2, 3, 1 ]	 \\
	\hline
	[ 2, 4, 3, 1, 5, 4, 3, 6, 5, 4, 2, 3, 4, 5 ]	 \\
	\hline
	[ 5, 4, 2, 3, 1, 6, 5, 4, 2, 3, 1, 4, 5, 6 ]	 \\
	\hline
	[ 2, 4, 5, 4, 2, 3, 1, 6, 5, 4, 2, 3, 1, 4 ]	 \\
	\hline
	[ 3, 1, 6, 5, 4, 2, 3, 1, 4, 3, 5, 4, 2, 6 ]	\\
\hline
\end{tabular}
\end{center}

}

\newpage
\paragraph{Table 8: Conjugacy Classes Associated with $e_6$ (continued)}

{
\footnotesize

\begin{center}
\begin{tabular}{|l|}
\hline\phantom{\rule{.2pt}{10pt}}
Conjugacy Classes \\
	\hline
	[ 2, 4, 2, 3, 1, 4, 5, 4, 3, 6, 5, 4, 2, 3, 1 ]	\\
	\hline
	[ 3, 1, 4, 3, 6, 5, 4, 2, 3, 1, 4, 3, 5, 4, 2 ]	\\
	\hline
	[ 3, 1, 4, 3, 5, 4, 2, 3, 1, 4, 3, 5, 4, 2, 6 ]	 \\
	\hline
	[ 2, 4, 3, 1, 5, 4, 2, 3, 4, 6, 5, 4, 2, 3, 4 ]	 \\
	\hline
	[ 2, 3, 1, 5, 4, 2, 3, 4, 5, 6, 5, 4, 2, 3, 4 ]	 \\
	\hline
	[ 2, 3, 1, 4, 2, 3, 1, 4, 5, 4, 2, 3, 4, 5, 6 ]	 \\
	\hline
	[ 2, 4, 2, 5, 6, 5, 4, 2, 3, 1, 4, 3, 5, 4, 2 ]	 \\
\hline
\end{tabular}
\quad
\begin{tabular}{|l|}
\hline\phantom{\rule{.2pt}{10pt}}
Conjugacy Classes \\
	\hline
	[ 2, 4, 3, 1, 5, 4, 2, 3, 1, 4, 3, 5, 6, 5, 4, 2 ]	 \\
	\hline
	[ 4, 2, 3, 1, 4, 3, 5, 4, 3, 1, 6, 5, 4, 2, 3, 1 ]	 \\
	\hline
	[ 2, 4, 2, 3, 1, 4, 3, 5, 4, 3, 1, 6, 5, 4, 2, 3 ]	 \\
	\hline
	[ 3, 1, 4, 2, 3, 4, 5, 4, 2, 3, 1, 4, 3, 6, 5, 4, 2 ]	 \\
	\hline
	[ 2, 3, 4, 3, 1, 6, 5, 4, 2, 3, 1, 4, 3, 5, 4, 2, 6 ]	 \\
	\hline
	[ 2, 5, 4, 3, 1, 6, 5, 4, 2, 3, 1, 4, 3, 5, 4, 2, 6 ]	\\
	\hline
	[ 2, 4, 3, 1, 5, 4, 2, 3, 4, 6, 5, 4, 2, 3, 1, 4, 3 ]	 \\
\hline
\end{tabular}
\end{center}

}

\newpage
\paragraph{Table 8: Conjugacy Classes Associated with $e_6$ (continued)}

{
\footnotesize

\begin{center}
\begin{tabular}{|l|}
\hline\phantom{\rule{.2pt}{10pt}}
Conjugacy Classes \\
	\hline
	[ 2, 3, 1, 4, 3, 5, 6, 5, 4, 2, 3, 1, 4, 3, 5, 4, 2, 6 ]	 \\
	\hline
	[ 2, 4, 5, 4, 2, 3, 1, 4, 3, 5, 4, 2, 6, 5, 4, 2, 3, 1 ]	 \\
	\hline
	[ 2, 3, 5, 4, 3, 1, 6, 5, 4, 2, 3, 1, 4, 3, 5, 4, 2, 6 ]	 \\
	\hline
	[ 2, 3, 1, 5, 4, 2, 3, 1, 4, 3, 6, 5, 4, 2, 3, 4, 5, 6 ]	 \\
	\hline
	[ 2, 3, 1, 5, 4, 3, 1, 6, 5, 4, 2, 3, 1, 4, 3, 5, 4, 6 ]	 \\
	\hline
	[ 2, 4, 3, 1, 5, 4, 6, 5, 4, 2, 3, 1, 4, 3, 5, 4, 2, 6 ]	\\
	\hline
	[ 2, 4, 3, 1, 5, 4, 2, 6, 5, 4, 2, 3, 1, 4, 3, 5, 4, 2 ]	\\
\hline
\end{tabular}
\end{center}

}

\newpage
\paragraph{Table 8: Conjugacy Classes Associated with $e_6$ (continued)}
{
\footnotesize

\begin{center}
\begin{tabular}{|l|}
\hline\phantom{\rule{.2pt}{10pt}}
Conjugacy Classes \\
	\hline
	[ 2, 4, 3, 1, 5, 4, 2, 3, 6, 5, 4, 2, 3, 1, 4, 3, 5, 4, 2, 6 ]
\\
	\hline
	[ 2, 3, 4, 2, 3, 1, 4, 3, 6, 5, 4, 2, 3, 1, 4, 3, 5, 4, 6, 5 ]
\\
	\hline
	[ 3, 1, 4, 2, 3, 4, 5, 6, 5, 4, 2, 3, 1, 4, 3, 5, 4, 2, 6, 5, 4 ]
\\
    \hline
	[ 3, 1, 5, 4, 2, 3, 1, 4, 3, 5, 6, 5, 4, 2, 3, 1, 4, 3, 5, 4, 2 ]
\\
	\hline
	[ 4, 3, 1, 5, 4, 2, 3, 4, 5, 6, 5, 4, 2, 3, 1, 4, 3, 5, 4, 2, 6 ]
\\
	\hline
	[ 2, 3, 1, 4, 2, 3, 1, 4, 5, 6, 5, 4, 2, 3, 1, 4, 3, 5, 4, 2, 6 ]
\\
	\hline
	[ 2, 3, 1, 4, 2, 3, 1, 4, 3, 5, 4, 2, 3, 1, 4, 3, 5, 6, 5, 4, 2 ]
\\
\hline

\end{tabular}
\end{center}

}

\newpage
\paragraph{Table 8: Conjugacy Classes Associated with $e_6$ (continued)}
{
\footnotesize

\begin{center}
\begin{tabular}{|l|}
\hline\phantom{\rule{.2pt}{10pt}}
Conjugacy Classes \\
	\hline
	[ 2, 3, 1, 4, 3, 5, 4, 3, 1, 6, 5, 4, 2, 3, 1, 4, 3, 5, 4, 2, 6, 5, 4 ]	 	
\\
	\hline
	[ 2, 4, 3, 1, 5, 4, 2, 3, 6, 5, 4, 2, 3, 1, 4, 3, 5, 4, 2, 6, 5, 4, 3 ]	 	
\\
	\hline
	[ 2, 3, 1, 4, 3, 5, 4, 2, 3, 6, 5, 4, 2, 3, 1, 4, 3, 5, 4, 2, 6, 5, 4 ]	 	
\\
	\hline
	[ 2, 3, 1, 4, 3, 5, 4, 2, 3, 1, 4, 3, 5, 4, 6, 5, 4, 2, 3, 1, 4, 5, 6 ]	 	
\\
	\hline
	[ 2, 3, 4, 2, 3, 1, 4, 6, 5, 4, 2, 3, 1, 4, 3, 5, 4, 2, 6, 5, 4, 3, 1 ]	 	
\\
	\hline
	[ 3, 1, 4, 3, 1, 5, 4, 2, 3, 4, 5, 6, 5, 4, 2, 3, 1, 4, 3, 5, 4, 2, 6 ]	 	
\\
	\hline
	[ 3, 1, 4, 3, 1, 5, 4, 2, 3, 1, 4, 5, 6, 5, 4, 2, 3, 1, 4, 3, 5, 4, 2 ]	 	
\\
\hline
\end{tabular}
\end{center}

}

\newpage
\paragraph{Table 9: Conjugacy Classes Associated with $e_7$}
{
\footnotesize

\begin{center}
\begin{tabular}{|l|}
\hline\phantom{\rule{.2pt}{10pt}}
Conjugacy Classes\\
	\hline
	[  ]	\\
	\hline
	[ 6 ]	\\
	\hline
	[ 2 ]	\\
	\hline
	[ 6, 5 ]\\
	\hline
	[ 4, 2 ]	\\
	\hline
	[ 2, 3 ]	\\
	\hline
	[ 2, 6 ]	\\
	\hline
	[ 4, 6 ]	\\
	\hline
	[ 4, 2, 3 ]	 \\
	\hline
	[ 5, 4, 2 ]	\\
	\hline
	[ 4, 2, 6 ]	 \\
	\hline
	[ 2, 6, 5 ]	 \\
	\hline
	[ 6, 5, 4 ]	 \\
	\hline
	[ 2, 3, 1 ]	 \\
	\hline
	[ 1, 4, 2 ]	 \\
	\hline
	[ 2, 3, 6 ]	 \\
	\hline
	[ 3, 1, 6 ]	 \\
	\hline
	[ 1, 4, 6 ]	 \\
	\hline
	[ 6, 5, 4, 2 ]\\
	\hline
	[ 4, 2, 3, 1 ]\\
	\hline
	[ 5, 4, 2, 3 ]\\
	\hline
	[ 1, 4, 2, 6 ]\\
	\hline
	[ 1, 5, 4, 2 ]\\
	\hline
	[ 4, 2, 3, 6 ]	 \\
	\hline
	[ 2, 3, 1, 6 ]\\
	\hline
	[ 3, 1, 6, 5 ]	 \\
	\hline
	[ 5, 4, 3, 1 ]	 \\
	\hline
	[ 4, 3, 1, 6 ]	 \\
	\hline
	[ 2, 3, 6, 5 ]	\\
	\hline
	[ 4, 2, 3, 1, 6 ]\\
	\hline
	[ 5, 4, 2, 3, 1 ]	\\
	\hline
	[ 1, 6, 5, 4, 2 ]	\\
	\hline
	[ 6, 5, 4, 2, 3 ]	 \\
	\hline
	[ 6, 5, 4, 3, 1 ]	\\
	\hline
	[ 2, 3, 1, 6, 5 ]	\\
	\hline
	[ 2, 5, 4, 2, 3 ]	\\
	\hline
	[ 2, 4, 3, 5, 4, 2 ]	\\
	\hline
	[ 2, 5, 4, 2, 3, 1 ]	\\
	\hline
	[ 6, 5, 4, 2, 3, 1 ]	\\
	\hline
	[ 4, 2, 5, 4, 2, 3 ]	\\
	\hline
	[ 2, 6, 5, 4, 2, 3 ]	\\
	\hline
	[ 2, 4, 3, 1, 5, 4, 2 ]	\\
	\hline
	[ 2, 4, 2, 3, 5, 4, 2 ]	\\
	\hline
	[ 2, 4, 3, 6, 5, 4, 2 ]	\\
	\hline
	[ 3, 1, 5, 4, 2, 3, 4 ]	\\
	\hline
    [ 2, 6, 5, 4, 2, 3, 1 ]	\\
	\hline
	[ 5, 4, 2, 6, 5, 4, 3 ]	\\
	\hline
	[ 2, 4, 3, 5, 4, 2, 3 ]	\\
	\hline
	[ 2, 4, 3, 1, 6, 5, 4, 2 ]	 \\
	\hline
	[ 2, 4, 2, 3, 5, 4, 2, 3 ]\\
	\hline
	[ 3, 1, 6, 5, 4, 2, 3, 4 ]	\\
	\hline
	[ 2, 3, 1, 4, 3, 5, 4, 2 ]	\\
	\hline
	[ 2, 4, 2, 3, 1, 5, 4, 2 ]	\\
	\hline
\end{tabular}
\quad
\begin{tabular}{|l|}
\hline\phantom{\rule{.2pt}{10pt}}
Conjugacy Classes\\
	\hline	
	[ 5, 4, 2, 6, 5, 4, 3, 1 ]	 \\
	\hline
	[ 2, 4, 2, 3, 6, 5, 4, 2 ]	 \\
	\hline
	[ 2, 4, 3, 6, 5, 4, 2, 3 ]	\\
	\hline
	[ 2, 3, 1, 4, 3, 6, 5, 4, 2 ]\\
	\hline
	[ 2, 4, 3, 1, 6, 5, 4, 2, 3 ]	 \\
	\hline
	[ 3, 1, 6, 5, 4, 2, 3, 1, 4 ]	 \\
	\hline
	[ 2, 4, 3, 1, 5, 4, 2, 3, 4 ]	 \\
\hline
\end{tabular}
\end{center}

}

\newpage
\paragraph{Table 9: Conjugacy Classes Associated with $e_7$ (continued)}
{
\footnotesize

\begin{center}
\begin{tabular}{|l|}
\hline\phantom{\rule{.2pt}{10pt}}
Conjugacy Classes \\
	\hline	
	[ 2, 3, 1, 5, 4, 2, 3, 1, 4, 6, 5, 4, 2, 3, 4 ]	 \\
	\hline
	[ 3, 1, 4, 3, 6, 5, 4, 2, 3, 1, 4, 3, 5, 4, 2 ]	\\
	\hline
	[ 2, 4, 2, 5, 6, 5, 4, 2, 3, 1, 4, 3, 5, 4, 2 ]	 \\
	\hline
	[ 2, 3, 4, 3, 1, 6, 5, 4, 2, 3, 1, 4, 3, 5, 4, 2 ]\\
	\hline
	[ 2, 4, 3, 1, 5, 4, 2, 3, 4, 5, 6, 5, 4, 2, 3, 4 ]	 \\
	\hline
	[ 2, 4, 3, 1, 5, 4, 2, 3, 4, 6, 5, 4, 2, 3, 1, 4 ]	 \\
	\hline
	[ 2, 4, 2, 3, 1, 6, 5, 4, 2, 3, 1, 4, 3, 5, 4, 2 ]	 \\
\hline
\end{tabular}
\end{center}

}

\newpage
\paragraph{Table 10: Conjugacy Classes Associated with $e_8$}
{
\footnotesize
\begin{center}
\begin{tabular}{|l|c|c|}
\hline\phantom{\rule{.2pt}{10pt}}
Conjugacy Classes C & Length & $|C \cap D_*(e)|$ \\
	\hline	
	[  ]	& 	0	& 	1 \\
	\hline
	[ 6 ]	& 	1	& 	36 \\
	\hline
	[ 6, 5 ]	& 	2	& 	240 \\
	\hline
	[ 4, 6 ]	& 	2	& 	1620 \\
	\hline
	[ 5, 4, 2 ]	& 	3	& 	270 \\
	\hline
	[ 4, 2, 6 ]	& 	3	& 	5184 \\
	\hline
	[ 2, 3, 6 ]	& 	3	& 	1440 \\
	\hline
	[ 6, 5, 4, 2 ]	& 	4	& 	4320 \\
	\hline
	[ 4, 3, 1, 6 ]	& 	4	& 	3240 \\
	\hline
	[ 3, 1, 6, 5 ]	& 	4	& 	480 \\
	\hline
	[ 5, 4, 2, 3 ]	& 	4	& 	540 \\
	\hline
	[ 2, 3, 1, 6 ]	& 	4	& 	1440 \\
	\hline
	[ 6, 5, 4, 3, 1 ]	& 	5	& 	45 \\
	\hline
	[ 4, 2, 3, 1, 6 ]	& 	5	& 	5184 \\
	\hline
	[ 2, 3, 1, 6, 5 ]	& 	5	& 	5760 \\
	\hline
	[ 5, 4, 2, 3, 1 ]	& 	5	& 	4320 \\
	\hline
	[ 6, 5, 4, 2, 3, 1 ]	& 	6	& 	720 \\
	\hline
	[ 4, 2, 5, 4, 2, 3 ]	& 	6	& 	1440 \\
	\hline
	[ 3, 1, 5, 4, 2, 3, 4 ]	& 	7	& 	4320 \\
	\hline
	[ 3, 1, 6, 5, 4, 2, 3, 4 ]	& 	8	& 	1440 \\
	\hline
	[ 2, 3, 4, 2, 3, 4, 5, 4, 2, 3, 4, 5 ]	& 	12	& 	6480 \\
	\hline
	[ 4, 2, 3, 1, 6, 5, 4, 2, 3, 1, 4, 5 ]	& 	12	& 	2160 \\
	\hline
	[ 2, 3, 1, 4, 2, 3, 4, 5, 4, 2, 3, 4, 5 ]	& 	13	& 	540 \\
	\hline
	[ 2, 3, 1, 4, 2, 3, 4, 6, 5, 4, 2, 3, 4, 5 ]	& 	14	& 	540 \\
	\hline
	[ 3, 1, 4, 3, 1, 5, 4, 2, 3, 1, 4, 5, 6, 5, 4, 2, 3, 1, 4, 3, 5, 4, 2, 6
 ]	& 	24	& 	80 \\
	\hline
\end{tabular}
\end{center}

}

\bigskip
%\newpage

\pagebreak
\vspace{5mm}

\noindent Zhenheng Li\\
Department of Mathematical Sciences \\
University of South Carolina Aiken\\
Aiken, SC 29801, USA\\
\noindent Email: zhenhengl@usca.edu\\

\noindent Zhuo Li\\
Department of Mathematics \\
Xiangtan University\\
Xiangtan, Hunan 411105, P. R. China\\
\noindent Email: zli@xtu.edu.cn\\

\end{document}